\documentclass[letterpaper, 10 pt, conference]{ieeeconf}  

\IEEEoverridecommandlockouts                              

\usepackage{amsmath}           
\usepackage{mathtools}
\usepackage{amssymb}
\usepackage{bbold}
\usepackage{ntheorem}
\usepackage[yyyymmdd,hhmmss]{datetime}

\usepackage{bm,upgreek}

\usepackage{graphics} 
\usepackage{epsfig} 
\usepackage{amsmath} 
\usepackage{amssymb}  
\usepackage{cite}
\usepackage{graphicx}
\usepackage{epstopdf}
\usepackage{color}

\newcommand{\Exp}[1]{\mathbb{E}[#1]}
\newcommand{\Prob}[1]{\mathbb{P}(#1)}          
\newcommand{\argmin}{\operatornamewithlimits{argmin}} 
\newtheorem{Assumption}{\textit{Assumption}}
\newtheorem{Lemma}{Lemma}
\theoremstyle{remark}
\newtheorem{Remark}{Remark}
\theoremstyle{Definition}
\newtheorem{Problem}{Problem}
\newtheorem{Theorem}{Theorem}
\newtheorem{Corollary}{Corollary}

\title{\LARGE \bf
Team Optimal Control of Coupled Subsystems with Mean-Field Sharing 
}

\author{Jalal Arabneydi and Aditya Mahajan
\thanks{This work was supported by the
Natural Sciences and Engineering Research Council of Canada through Grant NSERC-RGPIN 402753-11.}
\thanks{Jalal Arabneydi and Aditya Mahajan are with Department of Electrical Engineering, 
        McGill University, 3480 University St., Montreal, QC, Canada. 
        {\tt\small Email:jalal.arabneydi@mail.mcgill.ca} and
        {\tt\small Email:aditya.mahajan@mcgill.ca}}%
}

\begin{document}
\maketitle

\vspace*{-5.2cm}{\footnotesize{Proceedings of IEEE Conference on Decision and Control, 2014.}}
\vspace*{4.45cm}

\thispagestyle{empty}
\pagestyle{empty}

\begin{abstract}
We investigate team optimal control of stochastic subsystems that are weakly coupled in dynamics (through the mean-field of the system) and are arbitrary coupled in the cost. The controller of each subsystem observes its local state and the mean-field of the state of all subsystems. The system has a non-classical information structure. Exploiting the symmetry of the problem, we identify an information state and use that to obtain a dynamic programming decomposition. This dynamic program  determines a globally optimal strategy for all controllers. Our solution approach works for arbitrary number of controllers and generalizes to the setup when the mean-field is observed with noise.  The size of the information state is time-invariant; thus, the results generalize to the infinite-horizon control setups as well. In addition, when the mean-field is observed without noise, the size of the corresponding information state increases polynomially (rather than exponentially) with the number of controllers which allows us to solve problems with moderate number of controllers. We illustrate our approach by an example motivated by smart grids that consists of $100$ coupled subsystems.

\end{abstract}

\section{Introduction}

\subsection{Motivation}
Team optimal control of stochastic decentralized systems arises in many applications ranging from networked control systems, robotics,
communication networks, transportation networks, sensor networks, and economics. There is a long and rich history of research
on team theory, starting from the work of Radner \cite{Radner1962, MarschackRadner1972}, Witsenhausen \cite{Witsenhausen1968Counterexample,Witsenhausen1971separation, Witsenhausen1973}  and others; and continuing to various solution approaches that have been proposed in recent years. Due to space limitations, we can not provide a detailed overview of the literature; we rather refer the reader to \cite{Yuksel2013stochastic, Mahajan2012Information} for detailed overviews.

The scalability of the solution approach to large scale systems is an important consideration in team optimal control. Different approaches have been proposed to ensure that the solution complexity does not increase drastically with the number of subsystems. These include coordination-decomposition methods \cite{Culioli1990decomposition,Barty2010decomposition-coordination-heuristic} that use iterative message passing algorithm and mean-field games that reduce the optimal control problem to a game between an individual and the mass \cite{Caines2009springer,Gomes2013survey}, and references therein.

%
%
%

In this paper, we introduce a solution approach that exploits  symmetry to identify a low-dimensional information state. Our approach uses two steps. In the first step, we identify an equivalent  centralized system using the common information approach of \cite{Nayyar2013CIA}. In the second step, we exploit the symmetry of the system to identify an information state and use that to obtain a dynamic programming decomposition.


The rest of the paper is organized as follows. We formulate the team optimal control problem in Section \ref{Problem-Formulation} and identify the salient features of the problem and our contributions in Sections~\ref{Salient_features} and \ref{Contributions}, respectively. We present the main results in Section \ref{Main_results}, and provide some generalizations in Section \ref{Incomplete-Information}. In Section~\ref{Example}, we present an example (motivated by smart grid applications). 

%

\subsection{Notation}

To distinguish between random variables and their realizations, we use upper-case letters to denote random variables (e.g. $X$) and lower-case letters to denote their realizations (e.g. $x$). We use the short hand notation $X_{a:b}$  for the vector $(X_a,X_{a+1},\ldots,X_{b})$ and bold letters to denote vectors e.g. $\mathbf Y=(Y^1,\ldots,Y^n)$ where $n$ is the size of vector $\mathbf Y$.  $\mathbb{1}(\cdot)$ is the indicator function of a set, $\mathbb{P}(\cdot)$ is the probability of an event, $\mathbb{E}[\cdot]$ is the expectation of a random variable, and $|\cdot|$ is the cardinality of a set. $\mathbb{N}$ refers to the set of natural numbers.
\subsection{Problem Formulation}\label{Problem-Formulation}

Consider a discrete time decentralized control system with $n \in \mathbb{N}$ homogeneous subsystems that operate for a horizon $T \in \mathbb{N}$. The state of subsystem $i$, $i \in \{1,\ldots,n\}$, at time $t$, is denoted by $X^i_t \in \mathcal{X}$, where $\mathcal{
X}$ is a finite set (that does not depend on $i$). Let $U^i_t \in  \mathcal{U}$ denote the control action of controller  $i$, $i \in \{1,\ldots,n\}$, at time $t$, where $\mathcal{U}$ is a finite set (that does not depend on $i$).

We refer to the empirical distribution of all subsystems at time $t$ as the \textit{mean-field} of the system and denote it by $Z_t$, i.e., 
\begin{equation}
Z_t=\frac{1}{n} \sum_{i=1}^n \delta_{X^i_t}
\end{equation}
or equivalently, for $x \in \mathcal{X}$,
\begin{equation}\label{Mean-field}
Z_t(x)=\frac{1}{n}\sum_{i=1}^n \mathbb{1}(X^i_t=x)
\end{equation}
where $\delta_x$ denotes a Dirac measure on $\mathcal{X}$ with a point mass at $x$. Let $|\mathcal{X}|=k $ and $\mathcal{M}_n=\{(\frac{m_1}{n},\frac{m_2}{n},\ldots,\frac{m_k}{n}): m_i \in \{0,\ldots,n\}, \sum_{i=1}^k m_i=n\}$ denote the space of  realizations of $Z_t$. Note that $\mathcal{M}_n \subset \Delta(\mathcal{X})$, the space of probability distributions on $\mathcal{X}$.

\subsubsection{System dynamics}\label{System dynamics}

The subsystems are weakly coupled with each other in dynamics via the mean-field, as described below. The initial states of all subsystems are independent and distributed according to PMF (probability mass function) $P_X$ (that does not depend on $i$). The state $X^i_t$ of subsystem $i$ evolves according to 
\begin{equation}\label{Model}
X^i_{t+1}=f_t(X^i_t,U^i_t,W^i_t,Z_t), \hspace{.2cm} i \in \{1,\ldots,n\}
\end{equation}
where $f_t$ is the plant function at time $t$ and $\{W^i_t\}_{t=1}^T$ is an independent process with probability distribution $P_{W_t}$ at time $t$. Note that the plant functions $\{f_t\}_{t=1}^T \hspace{-.1cm}$ and the PMFs $\{P_{W_t}\}_{t=1}^T$ do not depend on $i$. 

The \textit{primitive random variables} $(X^1_1,\ldots,X^n_1,\{W^1_t\}_{t=1}^T,$ $\ldots,\{W^n_t\}_{t=1}^T)$ are mutually independent and defined on a common probability space.

\subsubsection{Information structure} \label{Information-structure}

In addition to the local state of its subsystems, each controller observes the history of the mean-field. Thus, the data available at controller $i$, $i \in \{1,\ldots,n\}$, at time $t$ is
\begin{equation}\label{Information_structure}
I^i_t=\{Z_{1:t},X^i_t\}.
\end{equation}

We refer to this information structure as \textit{mean-field sharing}. In Section \ref{Incomplete-Information}, we consider a generalization of this information structure in which each controller observes a noisy version of the mean-field. We refer to that information structure as \textit{partially observed mean-field sharing}.

The control action at controller $i$ is chosen according to
\begin{equation}\label{Complete-U}
U^i_t=g^i_t(Z_{1:t},X^i_t).
\end{equation}
The function $g^i_t$ is called the \textit{control law} of controller $i$ at time $t$. In this paper, we restrict attention to identical control laws at all controllers. In particular:

\begin{Assumption} \label{Homogeneous}
At any time $t$, the control laws at all controllers are identical i.e. $g^i_t=g^j_t$ for any $i,j \in \{1,\ldots,n\}$. Therefore, we drop the superscripts and denote the control law at every controller at time $t$ as $g_t$.
\end{Assumption}

In view of Assumption \ref{Homogeneous}, we call the collection $\mathbf{g}=(g_1,\ldots,g_T)$ of control laws over time as the \textit{control strategy} of the system.

\subsubsection{Cost-structure} \label{Cost-structure}

The subsystems are arbitrary coupled through cost. At each time step, the system incurs a cost that depends on joint state $\mathbf X_t=(X^1_t,\ldots,X^n_t)$ and joint action $\mathbf U_t=(U^1_t,\ldots,U^n_t)$ that is given by 
\begin{equation*}\label{Total-cost}
 \ell_t(\mathbf{X}_t,\mathbf U_t).
\end{equation*}
The performance of any strategy $\mathbf{g}$ is quantified by the expected total cost
\begin{equation}\label{Social-Cost}
 J(\mathbf{g})=\mathbb{E}^{\mathbf{g}}\left[\sum_{t=1}^T  \ell_t(\mathbf{X}_t,\mathbf U_t) \right]
\end{equation}
where the expectation is with respect to a joint measure induced on all system variables by the choice of $\mathbf{g}$.

\subsubsection{Optimization problem}
 We are interested in the following optimization problem.
\begin{Problem}\label{MDP}
Given the information structure in \eqref{Complete-U}, the horizon $T$, the plant functions $\{f_t\}_{t=1}^T$, the cost functions $\{\ell_t\}_{t=1}^T$, the PMF $P_X$ on the initial states, and the PMFs $\{P_{W_t}\}_{t=1}^T$ on the plant disturbance, identify a control strategy $\mathbf{g}^\ast$ to minimize the total cost $J(\mathbf{g})$ given by \eqref{Social-Cost}.
\end{Problem}

The above model assumes that all subsystems have access to the mean-field of the system. In certain applications such as cellular communications and smart grids, a centralized authority (such as a base station in cellular communication and an independent service operator in smart grids) may measure the mean field and transmit it to all controllers. In other applications such as multi-robot teams, all controllers may compute the mean-field in a distributed manner using methods such as consensus-based algorithms \cite{Olfati2006belief,Bishop2014}.

We first investigate the model where the mean field is shared perfectly and develop a solution methodology for that model. In Section \ref{Incomplete-Information}, we extend the solution methodology to a more practical model in which a noisy estimate of the mean field is observed.

\subsection{Salient Features of the Model}\label{Salient_features}

Our key simplifying assumption is that all control laws are identical (Assumption \ref{Homogeneous}). In general, this assumption leads to a loss in performance, as is illustrated by the example below.

\textit{Example}: Consider a system with $n$ homogeneous subsystems with control horizon $T=2$. Let state space and action space be $\mathcal{X}=\mathcal{U}=\{1,2,\ldots,n\}$ and probability distribution of initial states be uniform  on $\mathcal{X}$. Suppose that the system dynamics are given by 
\begin{equation}
X^i_2=U^i_1, \quad i\in \{1,\ldots,n\}.
\end{equation}

Let $\ell_1(\mathbf x_1,\mathbf u_1)=0$ and $\ell_2(\mathbf x_2,\mathbf u_2)=K\cdot\mathbb{1}(z_2 \neq \{\frac{1}{n},\frac{1}{n},\ldots,\frac{1}{n}\})$ 
 where $K$ is a positive number. The asymmetric strategy $ \bar{ \bm g}=(\bar{g}^1_1,\ldots,\bar{g}^n_1)$, where $\bar{g}^i_1(z_1,x^i_1)=i$, has a cost $J( \bar{\bm g})=0$. Hence, $ \bar{\bm g}$ is optimal. On the other hand, under any  symmetric strategy, $\Prob{\mathbb{1}(Z_2 \neq \{\frac{1}{n},\frac{1}{n},\ldots,\frac{1}{n}\})}$ is positive. Hence, a symmetric strategy is not globally optimal. By increasing $K$, we can make symmetric strategies perform arbitrary bad as compared to asymmetric strategies. 

Although assuming identical control laws (Assumption \ref{Homogeneous}) leads to loss in performance, it is a standard assumption in the literature on large scale systems for reasons of simplicity, fairness, and robustness. For example, similar assumption has been made in \cite{Schoute1978}, \cite{Shi2012survey}, \cite{Antsaklis2010}.

In the model described above, we assume that the strategies are pure (non-randomized). In general, randomized strategies are not considered in team problems because randomization does not improve performance \cite[Theorem 1.6]{Gihman1979}. However, if attention is restricted to identical strategies, randomized strategies may perform better than pure strategies \cite[Theorem 2.3]{Schoute1978}. In the above model, we assume that the control strategies are pure, primarily for the ease of exposition. As explained in the conclusion, our solution methodology generalizes to randomized strategies as well.

\subsection{Contributions } \label{Contributions}

In spite of the simplification provided by Assumption \ref{Homogeneous}, Problem~\ref{MDP} is conceptually challenging because it has a non-classical information structure \cite{Witsenhausen1971separation}. In general, team optimal control problems with non-classical information structure belong to NEXP complexity class \cite{Bernstein2002complexity}.  Although it is possible to get a dynamic programming decomposition for problems with non-classical information structure \cite{Witsenhausen1973}, the size of the corresponding information state increases with time. For some information structures, we can find information states that do not increase with time \cite{Mahajan2012Information}, but even for these models the size of the information state increases exponentially with the number of controllers.

Our key contributions in this paper are the following:
\begin{enumerate}
\item{We identify a dynamic program to obtain globally optimum control strategies.} 

\item{The size of the corresponding information state does not increase with time. Thus, our results extend naturally to infinite horizon setups.}

\item{The size of the corresponding information state increases  polynomially with the number of controllers. This allows us to solve problems with moderate number of controllers. (In Section \ref{Example}, we give an example with $n=100$ controllers).}

\item{The solution methodology and dynamic programming decomposition extend to the scenario where all controllers observe a noisy version of the mean-field.}

\end{enumerate} 


\section{ Main Results}\label{Main_results}

In this section, we use the common information approach \cite{Nayyar2013CIA} to introduce an equivalent  centralized problem (Problem~\ref{CIAMDP}) for Problem~\ref{MDP}. Then, we find an optimal solution for the equivalent problem and translate the obtained solution back to the solution of Problem~\ref{MDP}.

Following \cite{Nayyar2013CIA}, split the information $I^i_t$ available to controller $i$  into two parts: the \textit{common information} consisting of the history $Z_{1:t}$ of the mean-field process that is observed by all controllers; and the \textit{local information} consisting of the current state $X^i_t$ of subsystem $i$. Since the size of the local information does not increase with time, the model described above has a partial history sharing information structure \cite{Nayyar2013CIA}. For such systems, the structure of optimal control strategies and a dynamic programming decomposition was proposed in \cite{Nayyar2013CIA}. If we directly use these results on our model, the information state will be a posterior distribution on the global state $\mathbf{X}_t=(X^1_t,\ldots,X^n_t)$ of the system. As such the complexity of the solution increases doubly exponentially with the number of controllers.

To circumvent this issue, we proceed as follows.

\textit{Step 1:} We follow the common information approach proposed in \cite{Nayyar2013CIA} to convert the decentralized control problem into a centralized control problem from the point of view of a controller that observes the common information $Z_{1:t}$.

\textit{Step 2:} We exploit the symmetry of the problem (with respect to the controllers) to show that the mean-field $Z_t$ is an information state for the centralized problem identified in Step 1. We then use this information state $Z_t$ to obtain a dynamic programming decomposition.

The details of each of these steps are presented below.

\subsection{Step 1: An Equivalent Centralized System}\label{CIA}

Following \cite{Nayyar2013CIA}, we construct a fictitious centralized \textit{coordinated system} as follows. We refer to decision maker  in the coordinated system as the \textit{coordinator}. At time $t$, the coordinator observes the mean-field $Z_t$ and chooses a mapping  $\Gamma_t: \mathcal{X} \rightarrow \mathcal{U}$ as follows
\begin{equation}\label{Gamma-Psi}
\Gamma_t =\psi_t(Z_{1:t}).
\end{equation}

The function $\psi_t$ is called the \textit{coordination rule} at time $t$. The collection $\bm \psi=(\psi_1,\ldots,\psi_T)$ is called the \textit{coordination strategy}.

After the mapping $\Gamma_t$ is chosen, it is communicated to all controllers. Each controller in the coordinated system is a passive agent that uses its local state $X^i_t$ and the mapping $\Gamma_t$ to generate  
 \begin{equation}\label{Gamma-U}
U^i_t =\Gamma_t(X^i_t), \hspace{.2cm}  i\in \{1,\ldots,n\}.
\end{equation}
The dynamics of each subsystem and the cost function are the same as in the original problem. By a slight abuse of notation,  define
\begin{equation}
\ell_t(\mathbf{X}_t, \Gamma_t):= \ell_t(\mathbf X_t,\Gamma_t(X^1_t),\ldots,\Gamma_t(X^n_t)).
\end{equation}
The performance of any coordination strategy is quantified by the total expected cost 
\begin{equation}\label{JCIAMDP}
 \hat{J}(\bm {\psi})= \mathbb{E}^{\bm{\psi}}[\sum_{t=1}^T \ell_t(\mathbf X_t,\Gamma_t)]
\end{equation} 
where the expectation is with respect to a joint measure induced on all system variables by the choice of    
$\bm \psi$. 

Consider the following optimization problem.
\begin{Problem}\label{CIAMDP}
Given the information structure in \eqref{Gamma-Psi}, the horizon $T$, the plant functions $\{f_t\}_{t=1}^T$, the cost functions $\{\ell_t\}_{t=1}^T$, the PMF $P_X$ on the initial states, and the PMFs $\{P_{W_t}\}_{t=1}^T$ on the plant disturbance, identify a control strategy $\bm{\psi}^\ast$ to minimize the total cost $\hat{J}(\bm{\psi})$ given by \eqref{JCIAMDP}.
\end{Problem}

\begin{Lemma}[\cite{Nayyar2013CIA}, Proposition 3]\label{Equivalence}
Problem~\ref{MDP} and Problem~\ref{CIAMDP} are equivalent.
\end{Lemma}

 In particular, for any control strategy $\bm g=(g_1,\ldots,g_T)$ in Problem~\ref{MDP}, define a coordination strategy $\bm \psi=(\psi_1,\ldots,\psi_T)$ in Problem~\ref{CIAMDP} by
\begin{equation}\label{Psi}
\psi_t(z_{1:t}):=g_t(z_{1:t},\cdot), \quad \forall z_{1:t}.
\end{equation}
Then, $J(\bm g)=\hat{J}(\bm \psi)$. Similarly for any coordination strategy $\bm \psi$ in Problem~\ref{CIAMDP}, define a control strategy $\bm g$ in Problem~\ref{MDP} by
\begin{equation*}
g_t(z_{1:t}, x_t):=\psi_t(z_{1:t})(x_t), \quad \forall z_{1:t},\forall x_t.
\end{equation*}
Then, $J(\bm g)=\hat{J}(\bm \psi)$.

\subsection{ Step 2: Identifying an Information State and Dynamic Program }\label{Complete-Information}

An important result in identifying an information state is the following:

\begin{Lemma}\label{PX}
For any choice  $\gamma_{1:t}$ of $\Gamma_{1:t}$, any realization $z_{1:t}$ of $Z_{1:t}$, and any $\mathbf x \in \mathcal{X}^n$,
\[\begin{aligned}
\Prob{\mathbf X_t \hspace{-.1cm} = \hspace{-.1cm}\mathbf x |Z_{1:t} \hspace{-.1cm} = \hspace{-.1cm} z_{1:t},\Gamma_{1:t} \hspace{-.1cm} = \hspace{-.1cm} \gamma_{1:t}}&=\Prob{\mathbf X_t \hspace{-.1cm} = \hspace{-.1cm} \mathbf x |Z_t \hspace{-.1cm}= \hspace{-.1cm}z_t}\\
\hspace{-.1cm} &= \hspace{-.1cm} \frac{1}{|H(z_t)|}\mathbb{1}(\mathbf x \in H(z_t))
\end{aligned}\]
where $H(z)\hspace{-.1cm} := \hspace{-.1cm} \{\mathbf x  \in \hspace{-.1cm} \mathcal{X}^n \hspace{-.1cm} :\frac{1}{n} \sum_{i=1}^n \delta_{x^i}=z\}$.
\end{Lemma}
\textit{Proof outline:} To prove the result, it is sufficient to show that $\Prob{\mathbf X_t= \mathbf x|Z_{1:t}=z_{1:t},\Gamma_{1:t}=\gamma_{1:t}}$ is indifferent to permutation of $\mathbf{x}$. The latter can be proved using the symmetry of the model and the control laws. \hfill $\blacksquare$

Using this result, we can show that 

\begin{Lemma}\label{Lemma2}
The expected per-step cost may be written as a function of $Z_t$ and $\Gamma_t$. In particular, there exits a function $\hat{\ell}_t$ (that does not depend on strategy $\bm \psi$) such that 
\[\Exp{\ell_t(\mathbf X_t,\Gamma_t)|Z_{1:t},\Gamma_{1:t}}=:\hat{\ell}_t(Z_t,\Gamma_t).\]
\end{Lemma}
\textit{Proof outline:} Consider
\begin{equation*}
\begin{aligned}
&\Exp{\ell_t(\mathbf X_t,\Gamma_t)|Z_{1:t}=z_{1:t},\Gamma_{1:t}=\gamma_{1:t}}\\
&=\sum_{\mathbf{x}}\ell_t(\mathbf x,\gamma_t)\Prob{\mathbf{X}_t=\mathbf x|Z_{1:t}=z_{1:t},\Gamma_{1:t}=\gamma_{1:t}}.
\end{aligned}
\end{equation*}
Substituting the result of Lemma \ref{PX}, and simplifying gives the result. \hfill $\blacksquare$

\begin{Lemma}\label{PZ}
For any choice $\gamma_{1:t}$ of $\Gamma_{1:t}$, any realization $z_{1:t}$ of $Z_{1:t}$, and any $z \in \mathcal{M}_n$,
\[
\Prob{Z_{t+1} \hspace{-.1cm}= \hspace{-.1cm} z|Z_{1:t} \hspace{-.1cm} = \hspace{-.1cm} z_{1:t},\Gamma_{1:t} \hspace{-.1cm} =\hspace{-.1cm} \gamma_{1:t}} \hspace{-.1cm}= \hspace{-.1cm} \Prob{Z_{t+1} \hspace{-.1cm} =z|Z_t \hspace{-.1cm} = \hspace{-.1cm} z_{t}, \Gamma_t \hspace{-.1cm} =\hspace{-.1cm} \gamma_{t}}.\]
Also, above conditional probability does not depend on strategy $\bm \psi$.
\end{Lemma}
\textit{Proof outline:}  The result relies on the independence of the noise processes across subsystems and Lemma \ref{PX}. \hfill $\blacksquare$

Based on the results in steps 1 and 2, we have that


\begin{Theorem}\label{Thm-MDP}
In Problem~\ref{CIAMDP},  there is no loss of optimality in restricting attention to Markovian strategy i.e. $\Gamma_t=\psi_t(Z_t)$. Furthermore,  an optimal strategy $\bm \psi^\ast$  is obtained by  solving the following dynamic program. Define recursively value functions: 
\begin{equation}\label{DP-MDP1}
V_{T+1}(z_{T+1}):= 0, \quad \forall z_{T+1} \in \mathcal{M}_n
\end{equation}
and for $t= T,\ldots,1$, and for $ z_t \in \mathcal{M}_n$,  
\begin{equation}\label{DP-MDP2}
V_t(z_t):=\min_{\gamma_t}(\hat{\ell}_t(z_t,\gamma_t)+ \Exp{V_{t+1}(Z_{t+1})|Z_t=z_t,\Gamma_t=\gamma_t})
\end{equation}
where the minimization is over all functions $\gamma_t:\mathcal{X} \rightarrow \mathcal{U}$.  Let $\psi^\ast_t(z_t)$ denote any argmin of the right-hand side of  \eqref{DP-MDP2}. Then, the coordination strategy $\bm \psi^\ast=(\psi^\ast_1,\ldots,\psi^\ast_T)$ is optimal.

\end{Theorem}
\textit{Proof:}  $Z_t$ is an information state for Problem~\ref{CIAMDP} because:\\
1) As shown in Lemma~\ref{Lemma2}, the per-step cost can be written as a function of $Z_t$ and $\Gamma_t$.\\
2) As shown in  Lemma~\ref{PZ}, $\{Z_t\}_{t=1}^T$ is a controlled Markov process with control action $\Gamma_t$.\\
Thus, the result follows from standard results in Markov decision theory \cite{Bertsekas2012book}. \hfill $\blacksquare$

Based on the equivalence in Lemma~\ref{Equivalence}, we get

\begin{Corollary}
Let $\psi^\ast_t(z)$ be a minimizer of \eqref{DP-MDP2} at time $t$. Define 
\begin{equation}
g^\ast_t(z,x):=\psi^\ast_t(z)(x).
\end{equation}
Then, $\bm g^\ast \hspace{-.15cm}=\hspace{-.1cm}(g^\ast_1,\ldots,g^\ast_T)$ is an optimal strategy for Problem~\ref{MDP}.

%
%
\end{Corollary}

\begin{Remark}
\textnormal{The fictitious coordinated system is described only for ease of exposition. The dynamic program of \eqref{DP-MDP1} and \eqref{DP-MDP2} uses $z_t$ as the information state. Since $z_t$ is observed by each controller, each controller can independently solve the dynamic program; agreeing upon a deterministic rule to break ties while using $\argmin$ ensures that all controllers compute the same optimal strategy.}
\end{Remark}

\begin{Remark}
\textnormal{ The space $\mathcal{M}_n$ of realization of $z_t$ is finite and has cardinality less than $(n+1)^{|\mathcal{X}|}$. Thus, the solution complexity increases polynomially with the number of controllers.}
\end{Remark}
%

\section{Generalization To Partially Observed Mean-Field Sharing}\label{Incomplete-Information}

In this section, we consider a case where  mean-field is not completely observable. Let $Y_t \in \mathcal{Y}$ be a noisy measurement of $Z_t$ at time $t$ as follows:
\begin{equation}\label{Observation}
Y_t=h_t(Z_t,N_t)
\end{equation}
where $N_t$ is a random variable which takes value on a finite set $\mathcal{N}$. $\{N_t\}_{t=1}^T$ is an independent random process with PMF $P_{N_t}$, at time $t$, and is also mutually independent from all primitive random variables in  Section \ref{System dynamics}. Similar to \eqref{Complete-U}, we consider the following information structure:
\begin{equation}\label{Incomplete-U}
 U^i_t=g_t(Y_{1:t},X^i_t), \hspace{.2cm}  i\in \{1,\ldots,n\}
\end{equation}
 where $g_t:\mathcal{Y}^{t} \times \mathcal{X} \rightarrow \mathcal{U}$.

\begin{Problem}\label{POMDP}
Given the information structure in \eqref{Incomplete-U}, the horizon $T$, the plant functions $\{f_t\}_{t=1}^T$, the cost functions $\{\ell_t\}_{t=1}^T$, the PMF $P_X$ on the initial states, the PMFs $\{P_{N_t}\}_{t=1}^T$ on observation noise, and  the PMFs $\{P_{W_t}\}_{t=1}^T$ on the plant disturbance, identify a control strategy $\mathbf{g}^\ast$ to minimize the total cost $J(\mathbf{g})$ given by \eqref{Social-Cost}.
\end{Problem}

We follow the two-step approach of Section \ref{Main_results}. In step~1, we construct a centralized coordinated system in which a coordinator observes $Y_{1:t}$ and chooses
\begin{equation}\label{Gamma-Psi-incomplete}
\Gamma_t=\psi_t(Y_{1:t}).
\end{equation}
The rest of the setup is same as before. Similar to Problem~\ref{CIAMDP}, we get 

\begin{Problem}\label{CIAPOMDP} 
 Given the information structure in \eqref{Gamma-Psi-incomplete}, the horizon $T$, the plant functions $\{f_t\}_{t=1}^T$, the cost functions $\{\ell_t\}_{t=1}^T$, the PMF $P_X$ on the initial states, the PMFs $\{P_{N_t}\}_{t=1}^T$  on observation noise, and the PMFs $\{P_{W_t}\}_{t=1}^T$ on the plant disturbance, identify a control strategy $\bm{\psi}^\ast$ to minimize the total cost $\hat{J}(\bm{\psi})$ given by \eqref{JCIAMDP}.
\end{Problem}

As in Lemma~\ref{Equivalence}, Problem~\ref{POMDP} is equivalent  to Problem~\ref{CIAPOMDP}. In particular, for any control strategy $\bm g=(g_1,\ldots,g_T)$ in Problem~\ref{POMDP}, one can construct a coordination strategy $\bm \psi=(\psi_1,\ldots,\psi_T)$ in Problem~\ref{CIAPOMDP} that yields the same performance and vice versa.

In step 2, we show that $\Pi_t(z):=\Prob{Z_{t}=z|Y_{1:t}, \Gamma_{1:t-1}}$ is an information state for Problem~\ref{CIAPOMDP}. In particular:

\begin{Lemma}\label{Cost-POMDP}
There exists a function $\tilde{\ell}_t$ (that does not depend on strategy $\bm \psi$) such that
\begin{equation}
 \Exp{\ell_t(\mathbf{X}_t,\Gamma_t)|Y_{1:t},\Gamma_{1:t}}=:\tilde{\ell}_t(\Pi_t,\Gamma_t).
 \end{equation}
\end{Lemma}

\begin{Lemma} \label{Update-POMDP}
There exists a function $\phi_t$ (that does not depend on strategy $\bm \psi$) such that 
\begin{equation}\label{Update-Function}
\Pi_{t+1}=\phi_t(\Pi_t,\Gamma_t,Y_{t+1}).
\end{equation}
\end{Lemma}
Proofs of Lemma \ref{Cost-POMDP} and Lemma \ref{Update-POMDP} are omitted due to lack of space.
Similar to Theorem \ref{Thm-MDP}, we have that


\begin{figure*}[t!]
\centering{
\vspace{-.2cm}
\hspace*{-.0cm}
\scalebox{.87}{
\includegraphics[trim=3.5cm 0cm 1cm 0cm, clip, width=20cm ]{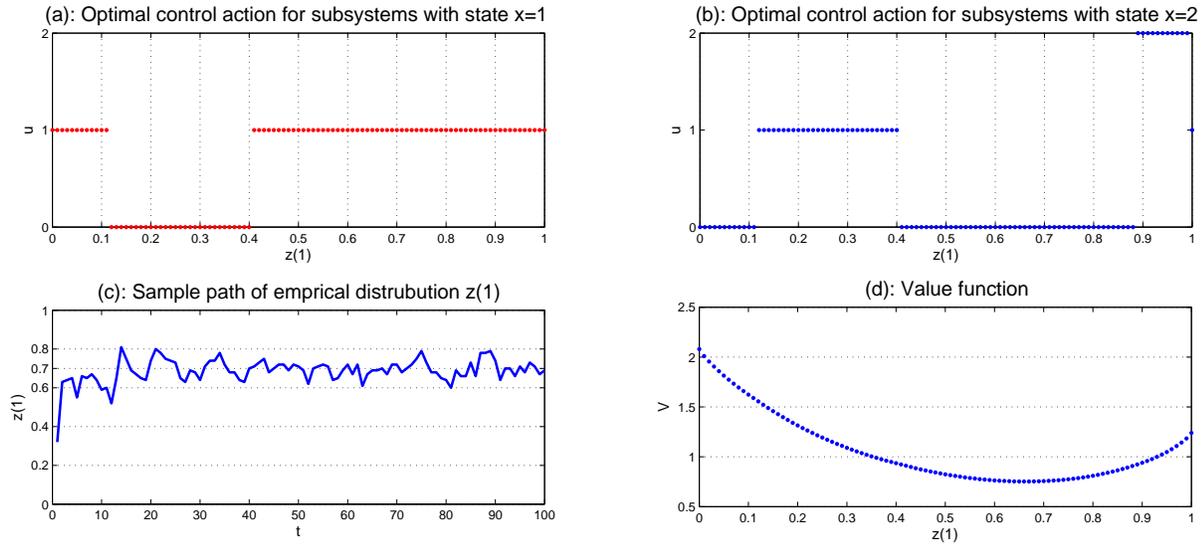}} 
\vspace{-.59cm}}
\caption{ Plots (a) and (b) show the optimal strategy as a function of $z(1)$. Plot (c)  shows the sample path of $z(1)$ for simulation time of $100$. Plot (d) depicts the value function with respect to $z(1)$. }
\label{fig_1}
\vspace{-.4cm}
\end{figure*}

\begin{Theorem}\label{Thm-POMDP}
 In  Problem~\ref{CIAPOMDP}, there is no loss of optimality in restricting attention to Markovian strategy i.e. $\Gamma_t=\psi_t(\Pi_t)$. Also,  optimal strategy $\bm \psi^\ast$ is obtained by  solving the following dynamic program. Let $\Delta(\mathcal{M}_n)$ denote the space of probability distributions on $\mathcal{M}_n$. Define recursively value functions:
\begin{equation}\label{DP-POMDP1}
V_{T+1}(\pi_{T+1})=0, \quad \forall \pi_{T+1} \in \Delta(\mathcal{M}_n)
\end{equation}
and for $t=T,\ldots,1$, and for $\pi_t \in \Delta(\mathcal{M}_n)$,
\begin{equation}\label{DP-POMDP2}
V_t(\pi_t)=\min_{\gamma_t}(\tilde{\ell}_t(\pi_t,\gamma_t)+ \Exp{V_{t+1}(\Pi_{t+1})|\Pi_t=\pi_t,\Gamma_t=\gamma_t})
\end{equation}
where the minimization is over all functions $\gamma_t:\mathcal{X} \rightarrow \mathcal{U}$. Let $\psi^\ast_t(\pi_t)$ denote any argmin of the right-hand side of  \eqref{DP-POMDP2}. Then, the coordination strategy $\bm \psi^\ast=(\psi^\ast_1,\ldots,\psi^\ast_T)$ is optimal.
\end{Theorem}
\textit{Proof:} $\Pi_t$ is an information state for Problem~\ref{CIAPOMDP} because:\\
1) As shown in Lemma~\ref{Cost-POMDP}, the expected per-step cost can be written as a function of $\Pi_t$ and $\Gamma_t$.\\
2) As shown in  Lemma~\ref{Update-POMDP}, $\{\Pi_t\}_{t=1}^T$ is a controlled Markov process with control action $\Gamma_t$.\\
Thus, the result follows from standard results in Markov decision theory \cite{Bertsekas2012book}. \hfill $\blacksquare$

Based on the equivalence between Problem~\ref{POMDP} and Problem \ref{CIAPOMDP}, we get


\begin{Corollary}
Let $\psi^\ast_t(\pi)$ be a minimizer of  \eqref{DP-POMDP2} at time $t$. Define 
\begin{equation}
g^\ast_t(\pi,x):=\psi^\ast_t(\pi)(x).
\end{equation}
Then, $\bm g^\ast \hspace{-.15cm}=\hspace{-.1cm}(g^\ast_1,\ldots,g^\ast_T)$ is an optimal strategy for Problem~\ref{POMDP}.

%
%
\end{Corollary}

%
%
%
%
%
%
%
%
%
%

\section{An Example}\label{Example}

In this section we consider an example of mean-field sharing that is motivated by applications in smart grids. Consider a system with $n$-devices where $\mathcal{X}=\{1,\ldots,k \}$ denotes the state space of each device and $\mathcal{U}=\{0,1,\ldots,k\}$ denotes the set of $k+1$ actions available at each device.

 Let $P(u)$ be  the controlled transition matrix under action $u \in \mathcal{U}$, i.e.
\begin{equation*}
[P(u)]_{xy}=\Prob{X^i_{t+1}=y \mid X^i_t=x, U^i_t=u},\quad x,y \in \mathcal{X}.
\end{equation*}

 Action $u=0$ is a \textit{free action} under which each device evolves in an uncontrolled manner, i.e. $P(0)=Q$, where $Q$ represents the \textit{natural} dynamics of the system. Action $u \neq 0$ is a \textit{forcing action} under which a fraction  $1-\epsilon_u$, $\epsilon_u \in [0,1]$, of devices switch to state $u$, and remaining $\epsilon_u$ devices follow the natural dynamics. Thus,

\begin{equation*}
P(u)=(1-\epsilon_u)\mathbf K_u + \epsilon_u Q
\end{equation*}
where $\mathbf K_u$ is a $k \times k$ matrix where column $u$ is all ones, and other columns are all zeros.

Action $u=0$ is free and it does not incur any cost, while action $u \neq 0$ incurs a cost $c(u)$. For notational convenience, let $c(0)=0$.

The objective is to keep the mean-field (i.e. the empirical distribution) of  the state of the devices close to a reference distribution $\zeta \in \Delta(\mathcal{X})$. The loss function is given by 
\begin{equation*}
\ell_t(\mathbf X_t, \mathbf U_t)=\frac{1}{n}\sum_{i=1}^n \hspace{-.1cm} c(U^i_t)+ D(Z_t \parallel  \zeta)
\end{equation*}
where $D(p \parallel q)$ denotes the  Kullback-Leibler divergence between $p,q \in \Delta(\mathcal{X})$ i.e. 
$D(p \parallel q)=\sum_{x \in \mathcal{X}} p(x) \log \frac{p(x)}{q(x)}.$  

The information structure is given by \eqref{Information_structure}. The objective is to choose a control strategy to minimize the infinite horizon discounted cost\footnote{Although we have only presented the details for finite horizon setup in this paper, the results generalize naturally to infinite horizon setup under standard assumptions. See Section \ref{Infinite_horizon} for a brief explanation.} 
\begin{equation}\label{Example_J}
J(\bm g)\hspace{-.1cm}= \mathbb{E}\hspace{-.1cm}\left[\sum_{t=1}^\infty \hspace{-.05cm}\beta^t \hspace{-.1cm}\left( \hspace{-.1cm} \frac{1}{n}\sum_{i=1}^n \hspace{-.1cm} c(U^i_t)+ D(Z_t \parallel  \zeta)\hspace{-.1cm}\right)\hspace{-.1cm}\right]
\end{equation}
where $\beta \in (0,1)$ is the discounted factor. 

A more elaborate variation of the above model is considered in \cite{Meyn2014ancillary} for controlling the operation of pool pumps.

Consider the above model for the following parameters
\begin{equation*}
\begin{aligned}
&n=100, \quad k=2, \quad \epsilon_1 =0.2, \quad \epsilon_2 =0.2,\\
& c(0)=0, \quad c(1)=0.1, \quad c(2)=0.2, \quad \beta=0.9,\\
&\zeta=\left[\begin{array}{c} 
 0.7 \\ 0.3 \end{array} \right],  \quad Q=\left[\begin{array}{cc} 
 0.25 \quad  0.75\\  0.375 \quad  0.625 \end{array} \right], \quad P_X=\left[\begin{array}{c} 
 \frac{1}{3} \\ \frac{2}{3} \end{array} \right].
\end{aligned}
\end{equation*}
The optimal time-homogeneous strategy for these parameters  is shown in Fig. \ref{fig_1}. Since state space is binary, $z(1)$ is sufficient to characterise the empirical distribution $z=[z(1), z(2)]$. Hence, for ease of presentation, we plot the optimal control law and value function as a function of the first component $z(1)$ of $z=[z(1), z(2)]$.

\section{Conclusion}
In this paper, we considered the team optimal control of decentralized systems with mean-field sharing. We follow a two-step approach: in the first step we construct an equivalent centralized system using the common information approach of \cite{Nayyar2013CIA}; in the second step, we exploit the symmetry of the system to identify information state and dynamic programming decomposition of the problem. We generalize our result to the case of partial observation of the mean field. We illustrate our results using an example motivated by smart grids. Our results extend naturally to the following setups.

\vspace*{-0.1cm}
\subsection{Randomized Strategies}\label{Randomized_strategies}

As mentioned earlier, if attention is restricted to identical control laws, then randomized strategies may perform better than pure strategies \cite[Theorem 2.3]{Schoute1978}. Our results extend naturally to randomized strategies by considering $\Delta(\mathcal{U})$, the space of probability distributions on $\mathcal{U}$, as the action space.
\vspace*{-0.1cm}

\subsection{Infinite Horizon}\label{Infinite_horizon}

The results of Lemma~\ref{Lemma2} and Lemma~\ref{PZ} are valid for the  infinite horizon setup as well. Hence, the results of Theorem~\ref{Thm-MDP} generalize to infinite horizon setup and under standard assumptions, the optimal coordination strategy is time-homogeneous and is given by the solution of a fixed point equation.
\vspace*{-0.1cm}

\subsection{Multiple Types of Subsystems}
We assumed that all subsystems are homogeneous. Consider a setup where subsystem $i$ has a type $\theta^i, \theta^i \in \Theta,$ and the dynamics are given by
$X^i_{t+1}=f_t(\theta^i,X^i_t,U^i_t,W^i_t,Z_t).$
Our results generalize to such a setup with $Z_t=\frac{1}{n} \sum_{i=1}^n \delta_{X^i_t,\theta^i}$.

\vspace*{-.1cm}

\bibliography{CDC_Ref}
\bibliographystyle{IEEEtran}

\newpage

\appendices

\section{ Proof of Lemma~\ref{PX}}\label{Proof_PX}

We use induction to prove the result. For notational convenience,  we denote $\Prob{A=a|B=b,C=c}$ by $\Prob{a|b,c}$.

Define $H(z):=\{\mathbf{x} \in \mathcal{X}^n: \frac{1}{n}\sum_{i=1}^n \delta_{x^i}=z\}$ as a set of all joint states $\mathbf x \in \mathcal{X}^n$ whose empirical distribution is $z$. Thus, at time $t$, we have 
\begin{equation}\label{Step1-0}
\mathbb{1}(z_t=\frac{1}{n}\sum_{i=1}^n \delta_{x^i_t})=\mathbb{1}(\mathbf{x}_t \in H(z_t)).
\end{equation} 
Notice that if $\mathbf{x}_t \in H(z_t)$, then one can interpret $H(z_t)$ as a collection of all permutations of $\mathbf x_t$ (such interpretation is critical for our proof).

In the first step, $t=1$, we have
\begin{equation}\label{Step1-1}
\begin{aligned}
&\Prob{\mathbf x_1|z_1,\gamma_1}\substack{{(a)}\\{=}}\Prob{\mathbf x_1|z_1}\substack{{(b)}\\{=}}\frac{\Prob{z_1|\mathbf x_1} \Prob{\mathbf x_1}}{\Prob{z_1}}\\
&=\frac{\mathbb{1}(z_1=\frac{1}{n}\sum_{i=1}^n \delta_{x_1^i})\Prob{\mathbf x_1}}{\Prob{z_1}}
\end{aligned}
\end{equation}
where $(a)$ follows from the fact that $\gamma_1=\psi_1(z_1)$ according to \eqref{Gamma-Psi} and $(b)$ follows from Bayes rule. From \eqref{Step1-0} and \eqref{Step1-1}, given $\{z_1,\gamma_1\}$, we get
\begin{equation}\label{Step1-7}
\Prob{\mathbf x_1|z_1,\gamma_1}=\begin{cases}
0 & \mathbf{x}_1 \notin H(z_1)\\
\alpha(z_1) & \mathbf{x}_1 \in H(z_1)
\end{cases}
\end{equation}
where $\alpha(z_1)=\frac{\Prob{\mathbf x_1}}{\mathbb{P}(z_1)}$ depends only on $z_1$. The reason lies in the fact that, when $\mathbf{x}_1 \in H(z_1)$, $H(z_1)$ contains nothing but permutations of $\mathbf x_1$ while joint probability distribution of initial states $\Prob{\mathbf{x}_1}=\prod_{i=1}^n P_X(x^i_1)$ is insensitive to permutation of $\mathbf{x}_1$. 
%
Since the summation of $\Prob{\mathbf x_1|z_1,\gamma_1}$ over $\mathbf x_1 \in \mathcal{X}^n$ is one, we have
\begin{equation}\label{Step1-5}
\alpha(z_1)=\frac{1}{|H(z_1)|}.
\end{equation} 
 From  \eqref{Step1-7} and \eqref{Step1-5}, we have
\begin{equation}
\Prob{\mathbf x_1|z_1,\gamma_1}=\Prob{\mathbf x_1|z_1}=\frac{\mathbb{1}(\mathbf x_1 \in H(z_1))}{|H(z_1)|}.
\end{equation}
Hence, the result holds for $t=1$. Assume the result holds for step $t$ i.e.
\begin{equation}\label{Step(t)-1}
\Prob{\mathbf x_t|z_{1:t},\gamma_{1:t}}=\Prob{\mathbf x_t|z_t}=\frac{\mathbb{1}(\mathbf x_t \in H(z_t))}{|H(z_t)|}.
\end{equation}
We prove that the result holds for step $t+1$ as follows.
\begin{equation}\label{Step(t+1)-1}
\begin{aligned}
&\Prob{\mathbf x_{t+1}|z_{1:t+1},\gamma_{1:t+1}}\substack{{(a)}\\{=}}\Prob{\mathbf x_{t+1}|z_{1:t+1},\gamma_{1:t}}\\
&\substack{{(b)}\\{=}}\frac{\Prob{z_{t+1}|\mathbf x_{t+1}} \Prob{\mathbf x_{t+1}|z_{1:t},\gamma_{1:t}}}{\Prob{z_{t+1}|z_{1:t},\gamma_{1:t}}}\\
&=\frac{\mathbb{1}(z_{t+1}=\frac{1}{n}\sum_{i=1}^n \delta_{x_{t+1}^i}) \Prob{\mathbf x_{t+1}|z_{1:t},\gamma_{1:t}}}{\Prob{z_{t+1}|z_{1:t},\gamma_{1:t}}}
\end{aligned}
\end{equation}
where $(a)$ follows from the fact that $\gamma_{t+1}=\psi_{t+1}(z_{1:t+1})$ according to \eqref{Gamma-Psi} and $(b)$ follows from Bayes rule. Similar to step $t=1$, we show that,  given $\{z_{1:t+1},\gamma_{1:t+1}\}$, above conditional probability is insensitive to permutation of $\mathbf{x}_{t+1}$. For that matter, we write the conditional probability in the numerator of \eqref{Step(t+1)-1} as follows.
\begin{equation}\label{Step(t+1)-2}
\begin{aligned}
&\Prob{\mathbf x_{t+1}|z_{1:t},\gamma_{1:t}}=\sum_{\mathbf x_t} \Prob{\mathbf x_{t+1},\mathbf x_{t}|z_{1:t},\gamma_{1:t}}\\
&=\sum_{\mathbf x_t} \Prob{\mathbf x_{t+1}|\mathbf x_{t},z_{1:t},\gamma_{1:t}}\Prob{\mathbf x_{t}|z_{1:t},\gamma_{1:t}}\\
&\substack{{(a)}\\{=}} \sum_{\mathbf x_t, \mathbf{w}_t}\left[\prod_{i=1}^n \mathbb{1}(x^i_{t+1}=f_t(x^i_t,\gamma_t(x^i_t),w^i_t,z_t)) \right]\\
&\cdot \Prob{\mathbf w_t} \cdot\Prob{\mathbf x_{t}|z_{1:t},\gamma_{1:t}}
\end{aligned}
\end{equation}
where $(a)$ follows from \eqref{Model} and the fact that $\mathbf{W}_t$ is independent from all data and decisions made before time $t$. Let $S:=\sigma(1,\ldots,n)$ denote an arbitrary permutation of set $\{1,\ldots,n\}$ and $S(i)$ denote the $i$th term of vector $S$. We use superscript $S$ to denote the permuted version of variables. For example, we denote $\mathbf x^S_{t+1}=(x_{t+1}^{S(1)},\ldots,x_{t+1}^{S(n)})$ as permuted version of $\mathbf{x}_{t+1}$ with respect to  $S$. Now, consider 
\begin{equation}\label{Step(t+1)-3}
\begin{aligned}
&\Prob{\mathbf x_{t+1}^S|z_{1:t},\gamma_{1:t}}=\hspace{-.2cm}\sum_{\mathbf x_t, \mathbf{w}_t} \hspace{-.1cm}\left[\prod_{i=1}^n \mathbb{1}(x^{S(i)}_{t+1}=f_t(x^i_t,\gamma_t(x^i_t),w^i_t,z_t)) \right]\\
&\cdot \Prob{\mathbf w_t} \cdot  \Prob{\mathbf x_{t}|z_{1:t},\gamma_{1:t}}\\
&\substack{{(a)}\\{=}}\sum_{\mathbf x_t, \mathbf{w}_t}\left[\prod_{i=1}^n \mathbb{1}(x^{S(i)}_{t+1}=f_t(x^{S(i)}_t,\gamma_t(x^{S(i)}_t),w^{S(i)}_t,z_t)) \right]\\
&\cdot \Prob{\mathbf w_t^S} \cdot  \Prob{\mathbf x_{t}^S|z_{1:t},\gamma_{1:t}},
\end{aligned}
\end{equation}
where $(a)$ follows from the fact that summation is insensitive to permutation. In particular, if  $D(\mathbf{x},\mathbf{w})$ is any arbitrary function of $(\mathbf{x}, \mathbf w)$, then we have
\begin{equation}
\sum_{\mathbf{x},\mathbf w} D(\mathbf x, \mathbf{w})=\sum_{\mathbf{x}^S, \mathbf{w}^S} D(\mathbf{x}^S,\mathbf{w}^S)=\sum_{\mathbf x, \mathbf w} D(\mathbf{x}^S,\mathbf{w}^S).
\end{equation}
Now, we consider terms in \eqref{Step(t+1)-3} separately as follows.

A) Since multiplication is insensitive to permutation,  the first term may be written as follows.
\begin{equation*} 
\prod_{i=1}^n \hspace{-.1cm} \mathbb{1}(\hspace{-.05cm}x^{\hspace{-.05cm}S(i)}_{t+1} \hspace{-.15cm}=\hspace{-.15cm} f_t\hspace{-.02cm}(\hspace{-.05cm}x^{\hspace{-.1cm}S(i)}_t\hspace{-.2cm},\gamma_t(\hspace{-.05cm}x^{\hspace{-.1cm}S(i)}_t\hspace{-.05cm}),\hspace{-.05cm}w^{\hspace{-.1cm}S(i)}_t\hspace{-.2cm},z_t\hspace{-.05cm})\hspace{-.05cm})\hspace{-.1cm}=\hspace{-.1cm}
\prod_{i=1}^n \hspace{-.1cm} \mathbb{1}(\hspace{-.05cm}x^i_{t+1}\hspace{-.15cm}=\hspace{-.15cm}f_t\hspace{-.02cm}(\hspace{-.05cm}x^i_t,\hspace{-.13cm}\gamma_t(\hspace{-.05cm}x^i_t\hspace{-.05cm}),\hspace{-.05cm}w^i_t,\hspace{-.05cm}z_t\hspace{-.05cm})\hspace{-.05cm})
\end{equation*} 

B) The second term may be written as follows. 
\begin{equation*}
\Prob{\mathbf w_t^S}=\prod_{i=1}^n P_{W_t}(w^{S(i)}_t)= \prod_{i=1}^n P_{W_t}(w^i_t)=\Prob{\mathbf w_t}.
\end{equation*}

C) According to \eqref{Step(t)-1}, the third term may be written as follows.
\begin{equation*}
\Prob{\mathbf x_{t}^S|z_{1:t},\gamma_{1:t}}\hspace{-.1cm}=\hspace{-.1cm}\frac{\mathbb{1}(\mathbf x_t^S \hspace{-.1cm}\in \hspace{-.1cm} H(z_t))}{|H(z_t)|}\hspace{-.1cm}=\hspace{-.1cm}\frac{\mathbb{1}(\mathbf x_t \hspace{-.1cm} \in \hspace{-.1cm} H(z_t))}{|H(z_t)|}\hspace{-.1cm}=\hspace{-.1cm}\Prob{\mathbf x_{t}|z_{1:t},\gamma_{1:t}}.
\end{equation*}

Substituting (A), (B), and (C) in \eqref{Step(t+1)-3}, we get
\begin{equation}\label{Step(t+1)-8}
\begin{aligned}
&\Prob{\mathbf x_{t+1}^S|z_{1:t},\gamma_{1:t}}\hspace{-.1cm}=\hspace{-.2cm} \sum_{\mathbf x_t, \mathbf{w}_t}\left[\prod_{i=1}^n \mathbb{1}(x^i_{t+1}\hspace{-.1cm}=\hspace{-.1cm}f_t(x^i_t,\gamma_t(x^i_t),w^i_t,z_t)) \right]\\
&\cdot \Prob{\mathbf w_t} \cdot  \Prob{\mathbf x_{t}|z_{1:t},\gamma_{1:t}}\substack{{(b)}\\{=}} \Prob{\mathbf x_{t+1}|z_{1:t},\gamma_{1:t}}
\end{aligned}
\end{equation}
where $(b)$ follows from \eqref{Step(t+1)-2}.

The rest of the proof is similar to that of  step $t=1$. From \eqref{Step(t+1)-1} and \eqref{Step(t+1)-8}, given $\{z_{1:t+1},\gamma_{1:t+1}\}$, we get
\begin{equation}\label{Step(t+1)-9}
\Prob{\mathbf x_{t+1}|z_{1:t+1},\gamma_{1:t+1}}=\begin{cases}
0 & \mathbf{x}_{t+1} \notin H(z_{t+1})\\
\alpha(z_{t+1}) & \mathbf{x}_{t+1} \in H(z_{t+1})
\end{cases}
\end{equation}
where $\alpha(z_{t+1})=\frac{\Prob{\mathbf x_{t+1}|z_{1:t},\gamma_{1:t}}}{\Prob{z_{t+1}|z_{1:t},\gamma_{1:t}}}$ depends only on $z_{t+1}$ because, when $\mathbf{x}_{t+1} \in H(z_{t+1})$, $H(z_{t+1})$ contains nothing but permutations of $\mathbf x_{t+1}$ while $\Prob{\mathbf x_{t+1}|z_{1:t},\gamma_{1:t}}$ is insensitive to permutation of $\mathbf{x}_{t+1}$ according to \eqref{Step(t+1)-8}. 
Since the summation of $\Prob{\mathbf x_{t+1}|z_{1:t+1},\gamma_{1:t+1}}$ over $\mathbf x_{t+1} \in \mathcal{X}^n$ is one, we have
\begin{equation}\label{Step(t+1)-7}
\alpha(z_{t+1})=\frac{1}{|H(z_{t+1})|}.
\end{equation} 
From  \eqref{Step(t+1)-9}  and \eqref{Step(t+1)-7}, we have
\begin{equation*}
\Prob{\mathbf x_{t+1}|z_{1:t+1},\gamma_{1:t+1}}\hspace{-.1cm}=\hspace{-.1cm}\Prob{\mathbf x_{t+1}|z_{t+1}}\hspace{-.1cm}=\hspace{-.1cm}\frac{\mathbb{1}(\hspace{-.05cm}\mathbf x_{t+1} \hspace{-.1cm}\in \hspace{-.1cm}H(z_{t+1}))}{|H(z_{t+1})|}.
\end{equation*}  \hfill $\blacksquare$
\vspace{-.1cm}
\section{Proof of Lemma~\ref{Lemma2}}\label{Proof_lemma2}
Consider the conditional expected cost at time $t$ given $\{z_{1:t},\gamma_{1:t}\}$ as follows
\[\begin{aligned}
&\Exp{\ell_t(\mathbf{X}_t,\Gamma_t)|Z_{1:t}=z_{1:t},\Gamma_{1:t}=\gamma_{1:t}}\\
&=\sum_{\mathbf x_t} \ell_t(\mathbf x_t, \gamma_t)\Prob{\mathbf X_t= \mathbf x_t| Z_{1:t}=z_{1:t},\Gamma_{1:t}=\gamma_{1:t}}\\
&\substack{{(a)}\\{=}}\sum_{\mathbf x_t} \ell_t(\mathbf x_t, \gamma_t)\frac{\mathbb{1}(\mathbf{x}_t \in H(z_t))}{|H(z_t)|}=:\hat{\ell}_t(z_t,\gamma_t)
\end{aligned}
\]
where $(a)$ follows from  Lemma \ref{PX}. Note that none of the above terms depend on strategy $\bm \psi$. \hfill $\blacksquare$\\

\vspace{-.4cm}

\section{Proof of Lemma~\ref{PZ}}\label{Proof_PZ}

For the sake of notational convenience, we denote $\Prob{A=a|B=b,C=c}$ by $\Prob{a|b,c}$ and $\mathcal{F}_t$ as the event $\{Z_{1:t}=z_{1:t}, \Gamma_{1:t}=\gamma_{1:t}\}$. Let $\alpha \in [0,1]$ and $x \in \mathcal{X}$, then we have
\begin{equation}\label{PZ1}
\allowdisplaybreaks
\begin{aligned}
&\Prob{Z_{t+1}(x)=\alpha |\mathcal{F}_t}=\mathbb{P}\left(\frac{1}{n} \left( \sum_{i=1}^n  \mathbb{1}(X^i_{t+1}=x)\right)=\alpha \middle\vert\ \mathcal{F}_t\right)\\
&\substack{(a)\\{=}}\sum_{\mathbf w,\mathbf x} \mathbb{1}\left(  \frac{1}{n} \left( \sum_{i=1}^n \mathbb{1}( f_t(x^i,\gamma_t(x^i), w^i, z_t)=x)\right)=\alpha \right)\\
&\cdot \Prob{\mathbf W_t=\mathbf w,\mathbf X_t=\mathbf x |\mathcal{F}_t}\\
&\substack{(b)\\=}\sum_{\mathbf w,\mathbf x} \mathbb{1}\left(  \frac{1}{n} \left( \sum_{i=1}^n \mathbb{1}( f_t(x^i,\gamma_t(x^i), w^i, z_t)=x) \right)=\alpha \right)\\
& \cdot \left[ \prod_{i=1}^n \Prob{W^i_t=w^i}\right]\Prob{\mathbf X_t=\mathbf x |\mathcal{F}_t}\\
&\substack{(c)\\=} \sum_{\mathbf w,\mathbf x} \mathbb{1}\left( \frac{1}{n} \left( \sum_{i=1}^n \mathbb{1}( f_t(x^i,\gamma_t(x^i), w^i, z_t)=x) \right)=\alpha \right)\\
& \cdot \left[ \prod_{i=1}^n \Prob{W^i_t=w^i}\right]\Prob{\mathbf X_t=\mathbf x |Z_t=z_t}\\
& =\Prob{Z_{t+1}(x)=\alpha|Z_t=z_t,\Gamma_t=\gamma_t}
\end{aligned}
\end{equation}
where $(a)$ follows from \eqref{Model} and the fact that ${\mathbf W_t}$ is independent from all data and decisions made before time $t$, and $(c)$ follows from Lemma~\ref{PX}. In addition, none the terms in \eqref{PZ1} depend on strategy $\bm \psi$.\hfill $\blacksquare$

\section{Proof of Lemma~\ref{Cost-POMDP}}\label{Proof_Cost_POMDP}
Consider the conditional expectation of per-step cost
\begin{equation*}\label{Conditional_cost_POMDP}
\begin{aligned}
&\Exp{\ell_t(\mathbf{X}_t,\Gamma_t)|Y_{1:t}=y_{1:t},\Gamma_{1:t}=\gamma_{1:t}}\\
&\substack{(a)\\=}\Exp{\hat{\ell}_t(Z_t,\Gamma_t)|Y_{1:t}=y_{1:t},\Gamma_{1:t}=\gamma_{1:t}}\\
&\substack{(b)\\=} \Exp{\hat{\ell}_t(Z_t,\Gamma_t)|Y_{1:t}=y_{1:t},\Gamma_{1:t}=\gamma_{1:t},\Pi_{1:t}=\pi_{1:t}}\\
&=\sum_{z_t} \hat{\ell}_t(z_t,\gamma_t)\Prob{Z_t \hspace{-.1cm}= \hspace{-.1cm} z_t|Y_{1:t} \hspace{-.1cm}= \hspace{-.1cm}y_{1:t},\Gamma_{1:t} \hspace{-.1cm} = \hspace{-.1cm}\gamma_{1:t},\Pi_{1:t} \hspace{-.1cm} = \hspace{-.1cm} \pi_{1:t}}\\
&=\sum_{z_t} \hat{\ell}_t(z_t,\gamma_t)\pi_t(z_t)=:\tilde{\ell}_t(\pi_t,\gamma_t)
\end{aligned}
\end{equation*}				
where $(a)$ follows from Lemma \ref{Lemma2} and $(b)$ follows from the fact that $\Pi_{1:t}$ is a function of $\{Y_{1:t},\Gamma_{1:t}\}$. Note that none of the above terms depend on strategy $\bm \psi$. \hfill $\blacksquare$

\section{Proof of Lemma~\ref{Update-POMDP}}\label{Proof_Update_POMDP}
For notational convenience, we denote $\Prob{A=a|B=b,C=c}$ by $\Prob{a|b,c}$. For $z_{t+1} \in \mathcal{M}_n$ and $y_{t+1} \in \mathcal{Y}$, we have
\begin{equation}\label{Lemma6.1}
\allowdisplaybreaks
\begin{aligned}
&\pi_{t+1}(z_{t+1})=\Prob{z_{t+1}|y_{1:t+1},\gamma_{1:t}}\substack{(a)\\=} \Prob{z_{t+1}|y_{1:t+1},\gamma_{1:t},\pi_{1:t}}\\
&=\frac{\Prob{z_{t+1},y_{t+1}|y_{1:t},\gamma_{1:t},\pi_{1:t}}}{\sum_{\tilde{z}_{t+1}} \Prob{\tilde{z}_{t+1},y_{t+1}|y_{1:t},\gamma_{1:t},\pi_{1:t}} }\\
&=\frac{\Prob{y_{t+1}|z_{t+1},y_{1:t},\gamma_{1:t},\pi_{1:t}}\Prob{z_{t+1}|y_{1:t},\gamma_{1:t},\pi_{1:t}}}{\sum_{\tilde{z}_{t+1}} \Prob{y_{t+1}|\tilde{z}_{t+1},y_{1:t},\gamma_{1:t},\pi_{1:t}}\Prob{\tilde{z}_{t+1}|y_{1:t},\gamma_{1:t},\pi_{1:t}}}\\
\end{aligned}
\end{equation}
where $(a)$ follows from the fact that $\Pi_{1:t}$ is a function of $\{Y_{1:t},\Gamma_{1:t}\}$. Consider the two terms of the denominator separately. The first term can be simplified as
\begin{equation}\label{Lemma6.3}
\begin{aligned}
&\mathbb{P}\hspace{-.05cm}(y_{t+1}\hspace{-.05cm}|\tilde{z}_{t+1},\hspace{-.05cm}y_{1:t},\hspace{-.1cm}\gamma_{1:t},\hspace{-.05cm} \pi_{1:t}\hspace{-.05cm})\\
&\substack{(b)\\=}\hspace{-.1cm}\sum_{n_{t+1}}\hspace{-.1cm}\mathbb{P}_{N_{t+1}}\hspace{-.05cm}(\hspace{-.05cm} n_{t+1} \hspace{-.05cm} )\hspace{-.0cm} \mathbb{1}(y_{t+1}\hspace{-.15cm}=\hspace{-.1cm}h_t(\tilde{z}_{t+1}, \hspace{-.05cm} n_{t+1} \hspace{-0.05cm})\hspace{-.05cm})=\Prob{y_{t+1}|\tilde z_{t+1}}
\end{aligned}
\end{equation}
where $(b)$ follows from \eqref{Observation}. The second term can be simplified as
\begin{equation}\label{Lemma6.2}
\begin{aligned}
&\Prob{\tilde{z}_{t+1}|y_{1:t},\gamma_{1:t},\pi_{1:t}}=\sum_{\tilde{z}_t}\Prob{\tilde{z}_{t+1},\tilde{z}_t|y_{1:t},\gamma_{1:t},\pi_{1:t}}\\
&=\sum_{\tilde{z}_t}\Prob{\tilde{z}_{t+1}|\tilde{z}_t,y_{1:t},\gamma_{1:t},\pi_{1:t}} \Prob{\tilde{z}_t|y_{1:t},\gamma_{1:t},\pi_{1:t}}\\
&\substack{(c)\\=}\sum_{\tilde{z}_t}\Prob{\tilde{z}_{t+1}|\tilde{z}_t,\gamma_t}\pi_t(\tilde{z}_t) =:\Prob{\tilde{z}_{t+1}|\pi_t,\gamma_t}
\end{aligned}
\end{equation}
where $(c)$ follows from Lemma~\ref{PZ} and definition of $\Pi_t$. From \eqref{Lemma6.1}, \eqref{Lemma6.3}, and \eqref{Lemma6.2}, we have
\begin{equation}\label{POMDP_proof_update}
\begin{aligned}
&\pi_{t+1}(z_{t+1})=\frac{\Prob{y_{t+1}|z_{t+1}}\Prob{z_{t+1}|\pi_t,\gamma_t}}{\sum_{\tilde{z}_{t+1}} \Prob{y_{t+1}|\tilde{z}_{t+1}}\Prob{\tilde{z}_{t+1}|\pi_t,\gamma_t}}\\
&=:\phi_t(z_{t+1},\pi_t, \gamma_t, y_{t+1})
\end{aligned}
\end{equation}
Furthermore,  none of the above terms  depend on strategy $\bm \psi$. \hfill $\blacksquare$

\end{document}